\documentclass[preprint, twocolumn, times, 5p]{elsarticle}

\usepackage[shortlabels]{enumitem}
\usepackage{amsmath,amssymb,amsfonts,mathtools}
\usepackage{dsfont}
\usepackage{color}
\usepackage{epstopdf}        
\usepackage{booktabs}
\usepackage{framed}
\usepackage{epsfig} 
\usepackage{hyperref}
\usepackage[dvipsnames]{xcolor}

\usepackage{diffcoeff}
\diffdef { p }
{
op-symbol = \partial ,
left-delim = \left .,
right-delim = \right | ,
subscr-nudge = +4 mu
}

\newtheorem{theorem}{Theorem}

\newtheorem{definition}{Definition}
\newtheorem{assumption}{Assumption}

\newenvironment{proof}{\paragraph{Proof}}{\hfill$\blacksquare$\vspace{5pt}}

\newcommand{\B}{\mathbb{B}}
\newcommand{\R}{\mathbb{R}}

\newcommand{\T}{\mathbb{T}}

\newcommand{\abs}[1]{\left\lvert#1\right\rvert}

\newcommand{\set}[1]{\left\{#1\right\}}

\newcommand*{\QEDB}{\hfill\ensuremath{\square}}%

\DeclareMathOperator*{\argmin}{arg\,min}

\allowdisplaybreaks

\usepackage{graphicx} 
\usepackage{graphics} 

\journal{~}

\begin{document}

\begin{frontmatter}



\title{\LARGE{\bf{High-Performance Optimal Incentive-Seeking in\\ Transactive Control for Traffic Congestion}}\tnoteref{t1,t2}}

\author{Daniel E. Ochoa \corref{cor1}}
\ead{daniel.ochoa@colorado.edu}

\author{Jorge I. Poveda}
\ead{jorge.poveda@colorado.edu}

\affiliation{organization={Department of Electrical, Energy and Computer Engineering. University of Colorado Boulder},
            city={Boulder},
            postcode={80305}, 
            state={Colorado},
            country={USA}}
            
\tnotetext[t1]{Research supported in part by NSF grant number CNS-1947613.}
\cortext[cor1]{Corresponding Author.}

\begin{abstract}
Traffic congestion has dire economic and social impacts in modern metropolitan areas. To address this problem, in this paper we introduce a novel type of model-free transactive controllers to manage vehicle traffic in highway networks for which precise mathematical models are not available. Specifically, we consider a highway system with managed lanes on which dynamic tolling mechanisms can be implemented in real-time using measurements from the roads. We present three incentive-seeking feedback controllers able to find in real-time the optimal economic incentives (e.g., tolls) that persuade highway users to follow a suitable driving behavior that minimizes a predefined performance index. The controllers are agnostic with respect to the exact model of the highway, and they are also able to guarantee fast convergence to the optimal tolls by leveraging non-smooth and hybrid dynamic mechanisms that combine continuous-time dynamics and discrete-time dynamics. We provide numerical examples to illustrate the advantages of the different presented techniques.
\end{abstract}


\end{frontmatter}

%
%
\section{INTRODUCTION}
The increase of population density in urban and sub-urban areas has triggered a significant growth of traffic congestion throughout the world, greatly affecting the commute of the public, as well as the expedited delivery of goods. 
For example, only in 2019, and solely in New York City, the economic losses induced by congestion climbed to \$11 billion USD. Moreover, commute times have significantly increased during the last years, forcing drivers to spend, on average, 41 hours per year in congested traffic during morning (6 am  to 9 am) and afternoon (3 pm to 6pm) peak travel times \cite{trafficJam}. This problem is only expected to  worsen during the next years, to the point that by the end of 2022 traffic congestion will cost $\$74$ billion USD to the economy of the United States. To tackle this challenge, cities throughout the world are developing and implementing automated control and optimization algorithms that can guarantee an optimal operation of the transportation infrastructure at all times. Examples include smart traffic light systems \cite{Kutadinata:14_Traffic}, dynamic pricing \cite{PhilipCSM,PovedaCDC17_a}, ride-sharing services \cite{Kleiner}, etc. Among these mechanisms, dynamic pricing has emerged as a promising technology to minimize congestion in dense cities such as London \cite{techreport_Tolling,london_case},  Milan \cite{MilanPricing}, and New York \cite{ny_pricing}. The goal of dynamic pricing is to induce ``optimal'' traffic flows that optimize a particular performance measure in the network by adaptively adjusting tolls or incentives \cite{pigou_book} based on the current state of the roads. To guarantee that the transportation system continuously operates at its optimal point, pricing algorithms must react \textsl{quickly} to changes in the traffic demand, weather conditions, road accidents, etc. This adaptability requirement has motivated the development of different recursive algorithms for optimal tolling computation, e.g., \cite{SandholmRouting2,TargetNashFlow,pricin_smart_grid,JohansoonTransporation}. Nevertheless, most existing pricing approaches are implemented based on (quasi) static lookup tables instead of real-time feedback traffic measurements, and therefore, under the presence of unexpected accidents or events in the system, are susceptible to generate sub-optimal or even ``perverse'' tolls that could exacerbate the very problems they were intended to solve \cite{PhilipCSM}.  Other recent approaches have relied on socio-technical models that aim to capture decision-making maps of drivers from recorded data; see \cite{annaswamy2018transactive}. In \cite{phan2016model}, \cite{zhang2013self} and \cite{zheng2016time}, the authors studied PID controllers to manage the operation of lanes in highway systems. In \cite{kachroo2016optimal} Hamilton-Jacobi-Bellman equations were solved for the optimal control of high-occupancy toll lanes, and adaptive algorithms based on linear parametrizations and welfare gradient dynamics were studied in \cite{poveda2017class}. A class of model-based saddle-flow dynamics were also recently studied in \cite{GB-JC-JP-ED:21-tcns} in the context of ramp metering control. For a recent review of transactive control for dynamic pricing see  \cite{lombardi2021model}.

In this paper, we depart from the traditional model-based approaches studied in the setting of transactive control, and instead, we introduce a new class of \emph{model-free} optimal incentive seeking controllers that can \emph{rapidly} learn optimal incentives (e.g., tolls) using only output measurements from the transportation systems, guaranteeing closed-loop stability at all times. Specifically, motivated by recent advances in non-smooth and hybrid extremum seeking control \cite{zero_order_poveda_Lina,PovedaKrsticFXT}, we introduce three incentive-seeking controllers (ISC) for model-free optimal price seeking in dynamic pricing: a smooth ISC that emulates the performance of a gradient-flow in the slowest time scale; a non-smooth ISC that emulates the behavior of fixed-time gradient flows in the slowest time scale; and a hybrid ISC that leverages momentum to improve transient performance in the slowest time scale. Each of the three controllers are interconnected with the dynamics of the highway, which incorporate socio-technical dynamics as well as traffic flows. Even though the controllers are agnostic to the traffic model, we establish practical asymptotic stability results for the resulting closed-loop system, and we numerically show that the non-smooth and hybrid ISCs can significantly outperform the smooth ISC in terms of transient performance under enough time scale separation in the closed-loop system. 

The rest of this paper is organized as follows: Section II presents preliminaries. Section III introduces the type of models that we consider in the paper. Section IV presents the proposed ISCs, as well as their stability results. Section V presents numerical examples and comparisons between all the controllers, and finally, Section VI ends with the conclusion.



%
\section{PRELIMINARIES}
\emph{Notation:} Given a compact set $\mathcal{A}\subset\R^n$ and a vector $z\in\R^n$, we use $|z|_{\mathcal{A}}\coloneqq \min_{s\in\mathcal{A}}\|z-s\|_2$ to denote the minimum distance of $z$ to $\mathcal{A}$. We use $\mathbb{S}^1\coloneqq \{z\in\R^2:z^2_1+z_2^2=1\}$ to denote the unit circle in $\R^2$, and $\mathbb{T}^n$ to denote the $n^{th}$ Cartesian product of $\mathbb{S}^1$. We also use $r\mathbb{B}$ to denote a closed ball in the Euclidean space, of radius $r>0$, and centered at the origin. We use $I_n\in\R^{n\times n}$ to denote the identity matrix, and $(x,y)$ for the concatenation of the vectors $x$ and $y$, i.e., $(x,y)\coloneqq [x^\top,y^\top]^\top$. A function $\beta:\R_{\geq0}\times\R_{\geq0}\to\R_{\geq0}$ is said to be of class $\mathcal{K}\mathcal{L}$ if it is non-decreasing in its first argument, non-increasing in its second argument, $\lim_{r\to0^+}\beta(r,s)=0$ for each $s\in\R_{\geq0}$, and  $\lim_{s\to\infty}\beta(r,s)=0$ for each $r\in\R_{\geq0}$.\\[0.1cm]
\emph{Hybrid Dynamical Systems: } In this paper, we will model our algorithms as Hybrid Dynamical Systems (HDS) with state $x\in \R^n$, and dynamics of the form:
\begin{subequations}\label{HDS}
\begin{align}
&x\in C,~\dot{x}= F(x),~~~~~x\in D,~x^+= G(x),
\end{align}
\end{subequations}
where $F:\mathbb{R}^n\to\mathbb{R}^m$ is called the flow map, and $G:\mathbb{R}^m\to\mathbb{R}^m$ is called the jump map. The sets $C$ and $D$, called the flow set and the jump set, define the points in $\mathbb{R}^n$ where the system can \emph{flow} or \emph{jump} according to $F$ or $G$, respectively. Thus, the HDS can be represented by the notation $\mathcal{H}=(C, F, D, G)$. Solutions  $x:\text{dom}(x)\to \R^n$ to \eqref{HDS} are parameterized by a continuous-time index $t\in\mathbb{R}_{\geq0}$, which increases continuously during flows, and a discrete-time index $j\in\mathbb{Z}_{\geq0}$, which increases by one during jumps. For a precise definition of hybrid time domains and solutions to \eqref{HDS} we refer the reader to \cite[Ch.2]{bookHDS}.\\[0.1cm]
\emph{Stability Notions:} The following definitions will be instrumental to characterize the convergence and stability properties of systems of the form \eqref{HDS}.
\begin{definition}
The compact set $\mathcal{A}\subset C\cup D$ is said to be \emph{uniformly asymptotically stable} (UAS) for system \eqref{HDS} if $\exists$ $\beta\in\mathcal{K}\mathcal{L}$ and $r>0$, such that for all solutions $x$ with $x(0,0)\in r\B$, and every $(t,j)\in \text{dom}(x)$, the following bound holds: $|x(t,j)|_{\mathcal{A}}\leq \beta(|x(0,0)|_{\mathcal{A}},t+j)$, $\forall~(t,j)\in\text{dom}(x)$. 
%
\QEDB
\end{definition}

We will also consider $\varepsilon$-parameterized HDS of the form: $x\in C_{\varepsilon}$, $\dot{x}=F_{\varepsilon}(x)$, $x\in D_{\varepsilon},~x^+\in G_{\varepsilon}(x)$, where $\varepsilon>0$. For these perturbed hybrid systems $\mathcal{H}_{\varepsilon}$ we will study \emph{practical stability} properties as $\varepsilon\to0^+$.
\begin{definition}\label{definitionPAS}
The compact set $\mathcal{A}\subset C\cup D$ is said to be \emph{practically Asymptotically Stable} (PAS) as $\varepsilon\to0^+$ for $\mathcal{H}_{\varepsilon}$ if $\exists$ $\beta\in\mathcal{K}\mathcal{L}$ and $r>0$ such that for each pair $(\tilde{\tau}, \nu)$ satisfying $r>\tilde{\tau}>\nu>0$, there exists $\varepsilon^*>0$ such that for all $\varepsilon\in(0,\varepsilon^*)$ every solution of $\mathcal{H}_{\varepsilon}$ with $|x(0,0)|_{\mathcal{A}}\leq \tilde{\tau}$ satisfies $|x(t,j)|_{\mathcal{A}}\leq \beta(|x(0,0)|_{\mathcal{A}},t+j)+\nu$, $\forall~(t,j)\in\text{dom}(x)$. 
~\hfill\QEDB
\end{definition}
The notions of P-AS can be extended to systems that depend on multiple parameters $\varepsilon=(\varepsilon_1,\varepsilon_2,\ldots,\varepsilon_{\ell})$. In this case, we say that $\mathcal{A}$ is P-AS as $(\varepsilon_{\ell},\ldots,\varepsilon_2,\varepsilon_{1})\to0^+$ where the parameters are tuned in order starting from $\varepsilon_1$.

Definitions 1 and 2 are similar to global and semi-global practical asymptotic stability properties studied in the literature \cite[Ch.7]{bookHDS}, with the difference that we restrict our attention only to initial conditions in a neighborhood of the set $\mathcal{A}$. These definitions are suitable for the application under study in this paper.
\begin{figure}[t]
    \centering
    \includegraphics[width=0.75\linewidth]{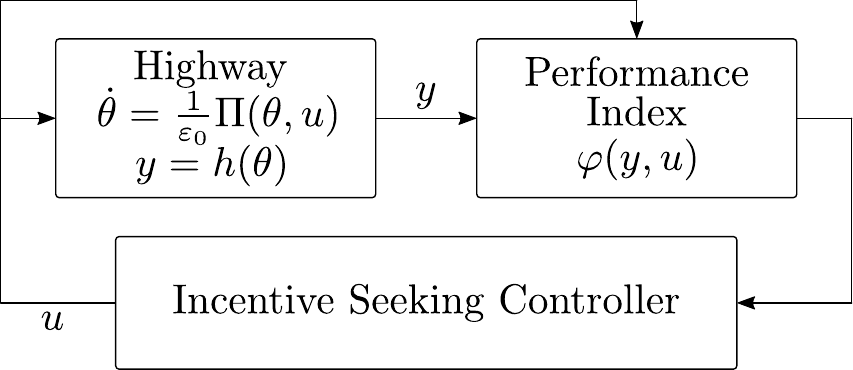}
    \caption{A closed-loop interconnection between an Incentive Seeking Controlller (ISC) and a highway system. The ISC will be designed to minimize in rel time a performance function defined by an external supervisor.}
    \label{fig:feedback}
\end{figure}
\emph{Problem Statement:} Consider a general highway network system modeled by a dynamical system of the form
\begin{equation}\label{ode:highway}
    \dot{\theta} = \frac{1}{\varepsilon_0} \Pi(\theta, u),~~~~~ y=h(\theta),\\
\end{equation}
where $\varepsilon_0$ determines the time scale of the dynamics, $\theta \in \R^n$ is the state of the highway, which can include the density of the cars per unit of length in a given lane, $u\in \R^m$ denotes exogenous incentives which can influence the behavior of highway users (e.g., tolls), and $y\in \R^p$ represents measurements that can be obtained from the highway state via the output map $h:\R^n\to \R^p$. Assume that an external supervisor or social planner provides a performance index $\varphi:\R^m\times\mathbb{R}^p\to\R$, which depends on the inputs and outputs of \eqref{ode:highway}. Our goal is to design feedback mechanisms able to find in real-time the optimal incentives that minimize the function $\varphi (y,u)$ at steady state. In particular, we consider closed-loop systems with the structure shown in Figure \ref{fig:feedback}, where the Incentive Seeking Controller (ISC) uses only real-time output \emph{measurements} of the performance index. The controller should be designed so that it can find the ``optimal'' incentives while preserving closed-loop stability at all times. In the following sections, we formalize each of the components illustrated in the scheme of Figure \ref{fig:feedback}.
\section{TRAFFIC IN HIGHWAY NETWORKS: SOCIO-TECHNICAL MODELS}\label{sec:model}
The performance of transportation systems is not solely dependent on their physical infrastructure, but also on their user behavior \cite{brown2017studies}. Indeed, in much of the literature that studies the modeling of dynamics in highway networks, the overall structure consists of a socio-technical model that combines a driver behavioral model and a traffic flow model; see \cite{lombardi2021model, engelson2006congestion, poveda2017class}. In this work, we follow a similar approach and we assume that the socio-technical and traffic flow models can be lumped together leading to highway network dynamics described by ODEs of the form \eqref{ode:highway}. Additionally, we make use of the following regularity assumption.
\begin{assumption}\label{highway:vc}
The map $\Pi(\cdot,\cdot)$ in \eqref{ode:highway} is locally Lipschitz. Moreover, there exists a compact set $\tilde{\Lambda}_\theta \coloneqq \lambda_\theta\B\subset \R^m$ with $\lambda_\theta \in \R_{>0}$, a closed set $\tilde{\Lambda}_u = \tilde{\Lambda}_u + \B$ where $\hat{\Lambda}_u\subset \R^m$, and a steady-state map $\ell:\R^m\to \R^n$ that is continuous and locally bounded relative to $\Lambda_u$, such that for each $\eta>0$ the compact set $\mathbb{M}_\eta \coloneqq \set{(\theta,u)~:~\theta =\ell(u),~u\in \Lambda_u\cap\eta\B,\theta\in\tilde{\Lambda}_{\theta}}$ is UAS for the HDS $\mathcal{H}_{ol} \coloneqq (\tilde{\Lambda}_\theta \times (\Lambda_u + \eta\B), \varepsilon_0^{-1}\Pi\times\set{0}, \emptyset, \emptyset)$ with state $(\theta, u)$.
\end{assumption}

In words, Assumption 1 guarantees that the highway dynamics are well-posed and stable with respect to external incentives $u$, and that the steady-state value of the traffic state is parameterized by $u$ via the map $\ell$. This assumption is standard (see \cite{poveda2017class}, \cite{JohansoonTransporation}, and \cite{SandholmRouting2}), and it is reasonable for many socio-technical models where external incentives $u$ determine the steady state equilibrium of the system.
\subsection{Socio-Technical Model}
To illustrate the advantages of the proposed ISC dynamics, we consider socio-technical models with a similar structure to the one described in \cite{annaswamy2018transactive}. In particular, we study the socio-technical model of a highway segment where drivers can choose between two parallel lanes: the general-purpose (GP) lane, which is uncharged, and the Express lane. Some of the motivations for choosing the Express lane include a faster travel time compared to the GP lane, as well as an expected reduced congestion. Prices (i.e., tolls) or subsidies can be assigned for the utilization of the Express lane depending on the traffic conditions. The model we consider focuses on the description of the average traffic density in the Express lane $\rho$, and the input flow of vehicles to the Express lane $q_{\text{EL}}$.  Figure \ref{fig:es_scheme0} shows a scheme representing the segment with the two parallel lanes:
\begin{figure}[t]
    \centering
    \includegraphics[width=0.9\linewidth]{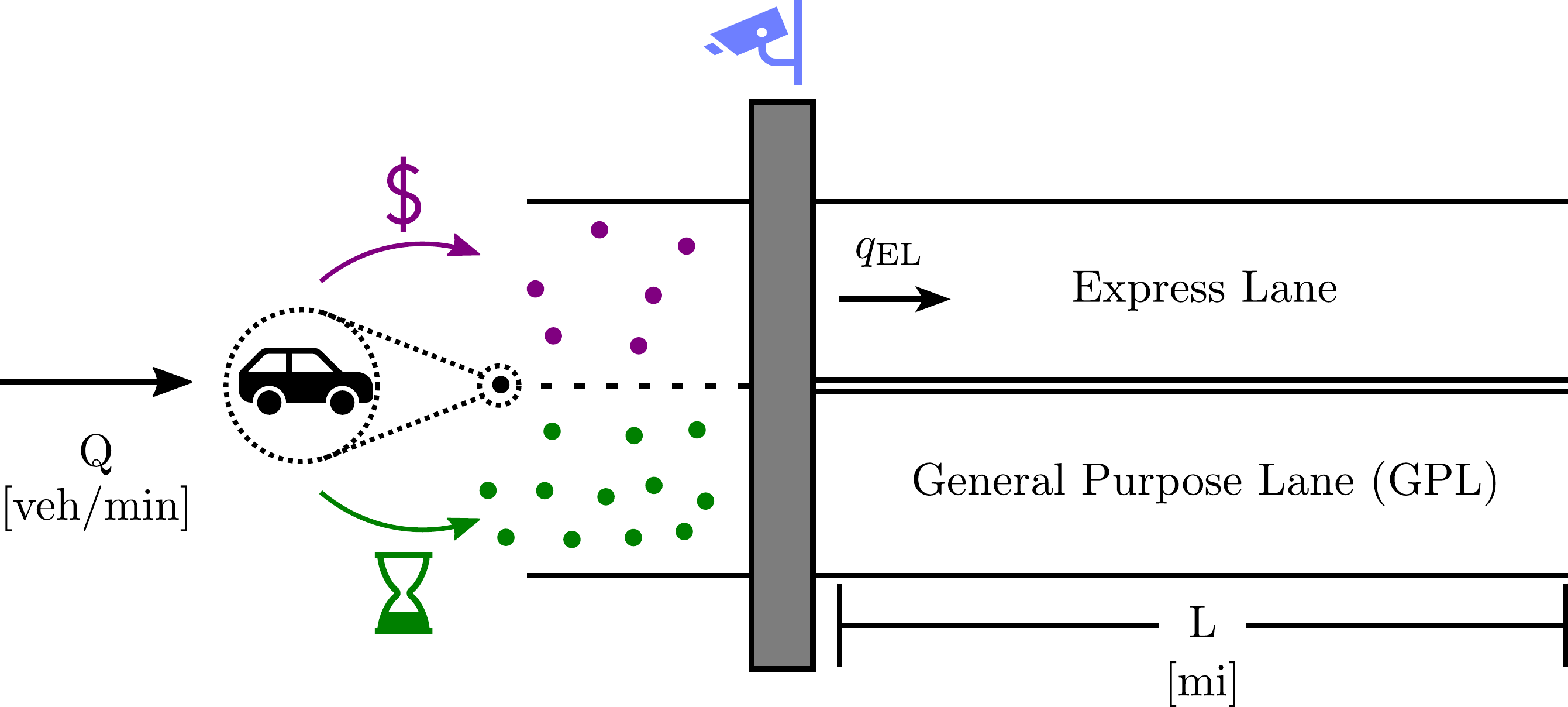}
    \caption{Scheme of segment with parallel lanes: the GP lane, and the Express lane.}
    \label{fig:es_scheme0}
\end{figure}
This model is divided into three main components:
\begin{enumerate}[a)]
    \item \emph{The driver behavioral model:} Each driver makes decisions based on the perceived costs of choosing either of the lanes. Some of the possible elements that can be included in these costs are travel times in the lanes, congestion, road conditions, among other quantities estimated by the drivers. To capture these scenarios, we model the costs by locally Lipschtiz functions $c_{\text{EL}}(q_{\text{EL}},\rho, u)$ and $c_{\text{GP}}(q_{\text{EL}},\rho)$ associated with the Express lane and the GP lane respectively,  where the input $u\in \R_{\ge 0}$ represents the tolls used to incentivize or to discourage the highway users from using the Express Lane.\\
    Naturally, the response of the individual drivers to the costs extends to the macroscopic level, concurrently affecting the input flow of vehicles to the Express Lane $q_{\text{EL}}$. In general, the drivers require a minimum transient time to adjust to changes in the marginal cost, which for instance could be induced by changes in the tolls. To account for this dynamic response, we model the dynamics of the macroscopic driver behavior as an ODE of the form
    \begin{equation}\label{sec:sociotechnical:ode}
        \dot{q}_{\text{EL}} = \Phi(q_{\text{EL}}, \rho, u),\quad q_{\text{EL}}\in[0,Q],
    \end{equation}
    where $\Phi:[0,Q]\times \R\times \R\to \R$ is a locally Lipschitz function that implicitly depends on the marginal cost of choosing the EL. With \eqref{sec:sociotechnical:ode}, we are able to capture a variety of social dynamics including, among others, evolutionary population dynamics whose stability properties have been recently studied in the literature, c.f. \cite{SandholmRouting2}.
    \item \emph{Equilibrium model:} The equilibrium model describes the resulting average velocity in the Express lane as a function of the average traffic density $\rho$. In this paper, we use a mollified version of the average velocity model presented in \cite{annaswamy2018transactive} and described by 
    \begin{equation}
        \overline{v}(\rho) = \frac{v_{\text{free}}-v_{\text{jam}}}{1+\text{exp}\left(\frac{4}{\rho_{\text{jam}}-\rho_{\text{critical}}}\left(\rho-\frac{\rho_{\text{jam}}+\rho_{\text{critical}}}{2}\right)\right)} + v_{\text{jam}},
    \end{equation}
    where $v_{\text{free}},v_{\text{jam}}$ are constants that represent the top speed and the jam vehicle speed in the Express lane, $\rho_{\text{critical}}$ denotes the average density below which the speed of the vehicles is expected to be close to $v_{\text{free}}$, and where $\rho_{\text{jam}}$ is the average vehicle density above which a traffic jam occurs in the Express lane. Consequently, these constants satisfy the relations: $v_{\text{free}}>v_{\text{jam}}$ and $\rho_{\text{critical}}<\rho_{\text{jam}}$.
    \item \emph{The traffic flow model:} This model represents the dynamics of the average traffic density $\rho$ of the Express lane, measured in vehicles per unit of length, as a function of the incoming rate of flow $q_{\text{EL}}$ and the average velocity of the Express lane $\overline{v}$. It is given by
    \begin{equation}
            \dot{\rho} = \frac{1}{L}\bigg(q_{\text{EL}} - \overline{v}(\rho)\rho\bigg),       
    \end{equation}
    where $L\in\R_{>0}$ represents the length of the highway segment under study. 
\end{enumerate}
By putting together the driver behavioral model and the traffic flow model, the dynamics of the average density in the express lane can be written in compact form as:
\begin{align}\label{dynamics:model}
    \dot{\theta}=\frac{1}{\varepsilon_0}\Pi(\theta, u){\coloneqq}\begin{pmatrix}
    k_{m}\Phi\left(q_{\text{EL}}, \rho, u\right)\\
    k_{\rho}\bigg(q_{\text{EL}}-\overline{v}(\rho)\rho\bigg)/L
\end{pmatrix},\quad y=h(\theta)
\end{align}
where $\theta\coloneqq (q_{\text{EL}}, \rho)$, and $h(\theta) \coloneqq \rho$. The ratio between the constants $k_m$ and $k_\rho$ in \eqref{dynamics:model}, dictates how fast the driver decisions occur in comparison with the overall traffic flow evolution described by $\rho$. In some cases, depending on the particular properties of the highway segment and the population of drivers, it might be the case that $k_m/k_\rho\gg 1$. For such scenarios, the relation between the driver response and the associated macroscopic behavior, captured by $q_{\text{EL}}$, can be simplified as a static map that depends on the marginal cost of choosing the Express lane:
    \begin{equation}
        q_{\text{EL}}(\rho, u) = \lambda(\tilde{c}_{\text{EL}}(\rho,u)-\tilde{c}_{\text{GP}}(\rho))Q,\label{sec:sociotechnical:qinstatic}
    \end{equation}
where $\lambda:\R\to [0,1]$ is a locally Lipschitz function that represents the traffic entering into the Express lane as a fraction of the total incoming traffic $Q$, which we measure in number of vehicles per amount of time, and where $\tilde{c}_{\text{EL}}$ and $\tilde{c}_{\text{GP}}$ are locally Lipschitz costs. When using relations of the form \eqref{sec:sociotechnical:qinstatic}, the socio-technical dynamics of \eqref{dynamics:model} is simplified as follows:
\begin{equation}\label{dynamics:model:qinstatic}
    \dot{\rho} = \frac{k_\rho}{L} \left(q_{\text{EL}}(\rho, u)-\overline{v}(\rho)\rho\right),\quad y = \rho.
\end{equation}
Note that \eqref{dynamics:model} and \eqref{dynamics:model:qinstatic} are particular cases of the ODE in \eqref{ode:highway}. For specific realizations of socio-technical models using static and dynamic formulations of the form \eqref{dynamics:model} and \eqref{dynamics:model:qinstatic}, we refer the reader to Section \ref{sec:numericalExamples}.
\subsection{Performance Indices}
\label{sec:performanceMeasures}
Depending on the objectives of the social planner, different performance indices can be considered for the purpose of real-time optimization. We will consider families of performance indices that satisfy the following assumption:
\begin{assumption}\label{assump:performanceIndex}
Suppose that Assumption \ref{highway:vc} holds, and let $\tilde{\varphi}(u):=\varphi(h(\ell(u)),u)$. The function $\tilde{\varphi}:\R^m\to \R$ is continuously differentiable, strictly convex in $\tilde{\Lambda}_u$, and its gradient is Lipschitz in $\tilde{\Lambda}_u$.
\end{assumption}

Sometimes, we will also use the following assumption:
\begin{assumption}\label{assumption:strongConvex}
There exists $\kappa>0$ such that $\tilde{\varphi}(\cdot)$ is $\kappa$-strongly convex in $\tilde{\Lambda}_u$ and its gradient is Lipschitz in $\tilde{\Lambda}_u$.
\QEDB
\end{assumption}

The above assumptions will guarantee enough regularity in the incentive-seeking problem, e.g., continuity of the cost and its gradient, the existence of finite optimal incentives, and sufficient monotonicity in the response map of the system. A particular example of a performance index satisfying Assumptions \ref{assump:performanceIndex} and \ref{assumption:strongConvex}, and that penalizes the deviation of the current vehicle-density $\rho$ from a desired operation point $\rho_{\text{ref}}$ provided by the external supervisor, is given by
\begin{equation}\label{reference:seekingcost}
 \varphi_{\text{ref}}(\theta, u)  = \abs{\rho-\rho_{\text{ref}}}^2.
\end{equation}
Performance indices of the form \eqref{reference:seekingcost} can be used by social planners who seek to improve traffic conditions, irrespective of the toll values needed to achieve such end. Among others, performance indices that explicitely depend on the toll prices $u$ can also be considered for situations in which the Express lane manager seeks profit maximization. In most of the cases, we will only require that Assumptions \ref{assump:performanceIndex} or \ref{assumption:strongConvex} hold in a neighborhood of the minimizer of $\tilde{\varphi}$. 
\begin{figure}[t]
    \centering
    \includegraphics[width=0.8\linewidth]{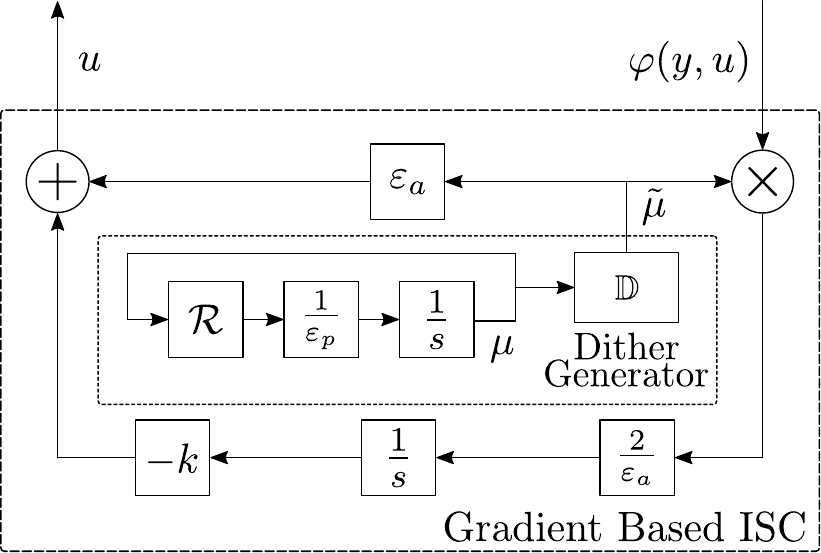}\vspace{-5pt}
    \caption{Gradient Based Incentive Seeking Control}
    \label{fig:es_scheme}
\end{figure}
\section{INCENTIVE SEEKING FEEDBACK SCHEMES}\label{sec:isf}
In this section, we introduce three different ISC algorithms able to guarantee convergence and stability of the set
$
\mathcal{A}_0:=\{(\theta, u)\in \tilde{\Lambda}_\theta\times \tilde{\Lambda}_u:\theta=\ell(u),u=\argmin_{u}\tilde{\varphi}(u)\},
$
where $\Lambda_u$ and $\tilde{\Lambda}_\theta$ are compact sets, $\ell$ is given in Assumption \ref{highway:vc}, and $\tilde{\varphi}$ is as in Assumption \ref{assump:performanceIndex}.
The ISCs make use of small exploration signals injected into the transportation dynamics for the purpose of real-time learning. These signals are generated by dynamic oscillators of the form:
\begin{equation}\label{dither}
    \dot{\mu} =\frac{1}{\varepsilon_p}\mathcal{R}\mu,\qquad\mu\in \T^m,
\end{equation}
where $\varepsilon_p\in \R_{>0}$ is a tunable parameter, and the matrix $\mathcal{R}\in \R^{2n\times 2n}$ is a block diagonal matrix with blocks given by $\mathcal{R}_i = 2\pi\begin{pmatrix} 0           &\omega_i\\-\omega_i   & 0\end{pmatrix}$, with $\omega_i\in \R_{>0}$ and for $i=\set{1,\cdots, n}$. We use $\omega\coloneqq(\omega_1,\cdots, \omega_n)$ to denote the vector of frequencies of the signals, and we consider ISCs that generate incentives $u$ of the form
%
%
\begin{equation}\label{inputAndmapMu}
    u= \hat{u} + \varepsilon_a \mathbb{D}\mu,
\end{equation}
where $\mathbb{D}\mu = (\mu_1, \mu_3, \cdots, \mu_{2n-1})$ represents the odd components of $\mu$, and $\hat{u}$ is the nominal incentive generated by each particular algorithm. We will impose the following assumption on $\omega$.
\begin{assumption}\label{assumption:ditherfreq}
The dithering frequencies $\omega_i$ satisfy: 1) $\omega_i>0$ is a rational number for all $i$, and 2) there are no repeated dither frequencies, i.e., $i\neq j \implies \omega_i\neq\omega_j$ and $\omega_i\neq 2\omega_j$. 
\end{assumption}
Assumption \eqref{assumption:ditherfreq} guarantees orthogonality conditions for the dither signals used by the ISCs to update the incentives $u$. These conditions will enable real-time learning in the closed-loop system via averaging theory.
\subsection{Gradient Based Incentive Seeking Control}
We first consider a smooth ISC, denoted GISC, presented in Figure \ref{fig:es_scheme}, which generates the nominal incentive $\hat{u}$ via the following differential equation:
\begin{equation}\label{gisf}
    \left(\begin{array}{c}
         \dot{\hat{u}} \\
         \dot{\mu}
    \end{array}\right)= F_1(x_1) \coloneqq \begin{pmatrix}
    -k\varphi (y,u)M(\mu)\\
    \frac{1}{\varepsilon_p}\mathcal{R}\mu
    \end{pmatrix},~ x_1\in \R^{m}\times \T^m,
\end{equation}
where $x_1\coloneqq (\hat{u}, \mu)$, and $M(\mu) = \frac{2}{\varepsilon_a}\mathbb{D}\mu$. The controller \eqref{gisf} is based on smooth extremum-seeking controllers \cite{KrsticBookESC}, which aim to emulate gradient flows whenever the highway dynamics \eqref{ode:highway} are neglected. The controller makes use of direct measurements of the perfomance index $\varphi(y,u)$, and therefore it is agnostic to the dynamics of the transportation system. In the context of traffic congestion, related dynamics have been studied in \cite{poveda2017class} for adaptive pricing in affine congestion games, \cite{krsticPDE} for highways with bottlenecks, and in \cite{sanchez2020dynamic} via simulations for congestion lanes.  The following theorem shows that \eqref{gisf} is a suitable controller to learn optimal incentives in transportation systems with socio-technical dynamics in the loop.  
\begin{theorem}\label{thm:gisf}
Suppose that Assumptions \ref{highway:vc}, \ref{assump:performanceIndex} (or \ref{assumption:strongConvex}) and \ref{assumption:ditherfreq} hold. Then, the closed-loop system corresponding to Figure \ref{fig:feedback} with ISC given by \eqref{gisf}, renders PAS the set $\mathcal{A}_1 \coloneqq \mathcal{A}_0\times \mathbb{T}^m$ as $(\varepsilon_0,\varepsilon_p,\varepsilon_a)\to 0^+$. 
\QEDB
\end{theorem}
\vspace{-0.4cm}
\begin{proof}
The result of Theorem \ref{thm:gisf} can be established by showing that all the assumptions needed to apply \cite[Thm.1]{poveda2017framework} are satisfied in a neighborhood of the optimal incentive. First, note that, by Assumption \ref{highway:vc}, the plant has a well-defined steady state input-to-output map $\tilde{\varphi}$. Also, under Assumption \ref{assump:performanceIndex}, this response map is strictly convex, and thus has a unique minimizer. Since under Assumption \ref{assumption:ditherfreq} the average dynamics of \eqref{gisf} can be computed to be $\dot{\hat{u}}=-k\nabla \tilde{\varphi}(u)+\mathcal{O}(\varepsilon_a)$ (see, e.g., \cite[Sec. 7]{zero_order_poveda_Lina}), it follows that for $\varepsilon_a$ sufficiently small, in a neighborhood of the optimal incentive $u^*$ the average dynamics converge to a neighborhood of $u^*$. By averaging theory and the results of \cite[Thm.1]{poveda2017framework}, the original system retains the stability properties in a practical sense. The result follows by using a (local) singular perturbation argument to interconnect the dynamics \eqref{gisf} with the dynamics \eqref{dynamics:model}.
\end{proof}
\begin{figure}[t]
    \centering
    \includegraphics[width=0.8\linewidth]{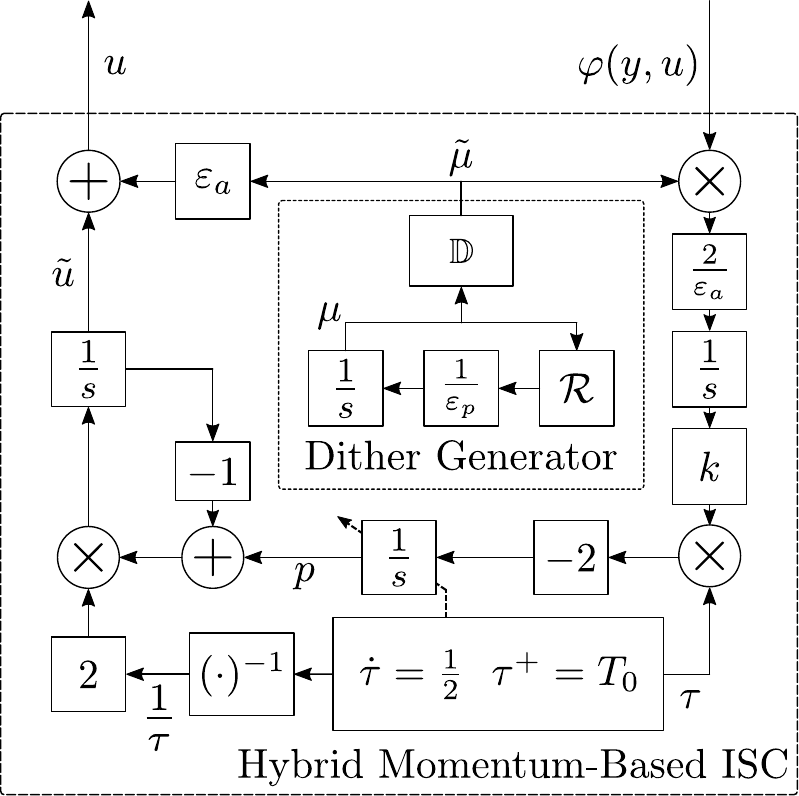}\vspace{-5pt}
    \caption{Hybrid Momentum Based Incentive Seeking Control}
    \label{fig:hmes_scheme}
\end{figure}
%

While the ISC \eqref{gisf} can achieve optimal incentive seeking, as $(\varepsilon_0,\varepsilon_p,\varepsilon_a)\to0^+$ the rate of convergence achieved by this controller emulates the convergence rate of a gradient descent flow, which is either of order $\mathcal{O}(1/t)$ or $\mathcal{O}(e^{-\kappa t})$, where $\kappa$ is given by Assumption \ref{assumption:strongConvex} (note that Assumption \ref{assumption:strongConvex} implies Assumption \ref{assump:performanceIndex}). In Section 5, we will show that the steady-state performance function related to the socio-technical model of the traffic network can have drastically different ``slopes'' near the optimal point, including cases where $\kappa\ll1$. These ``flat'' regions can drastically deteriorate the performance of controllers that seek to emulate traditional gradient flows, e.g., system \eqref{gisf}. To achieve better transient performance in this scenario, we now consider a class of hybrid dynamics that use momentum.
\subsection{Hybrid Momentum-Based Incentive Seeking Control}
To achieve better transient performance compared to \eqref{gisf}, we now consider the hybrid ISC shown in Figure \ref{fig:hmes_scheme}, which has continuous-time and discrete-time dynamics given by:
\begin{subequations}\label{hmisf}
\begin{align}
    &\left(\begin{array}{c}
         \dot{\hat{u}} \\
         \dot{p}\\
         \dot{\tau}\\
         \dot{\mu}
    \end{array}\right)= F_2(x_2) \coloneqq \begin{pmatrix}
    \frac{2}{\tau}\left(p - \hat{u}\right)\\
    -2k\tau \varphi (y,u)M(\mu)\\
    \frac{1}{2}\\
    \frac{1}{\varepsilon_p}\mathcal{R}\mu
    \end{pmatrix},\label{himsf:flowMap}\\
    &x_2 \in  C_2\coloneqq \Big\{x_2\in \R^{2m + 1}\times \T^m~:~\tau \in [T_0,T ]\Big\},\label{himsf:flowSet}\\
    &\left(\begin{array}{c}
         \hat{u}^+ \\
         p^+\\
         \tau^+\\
         \mu^+
    \end{array}\right)= G_2(x_2) \coloneqq \begin{pmatrix}
    \hat{u}\\
    \sigma p + (1-\sigma)q\\
    T_0\\
    \mu
    \end{pmatrix},\label{himsf:jumpMap}\\
    &x_2 \in D_2\coloneqq \Big\{x\in \R^{2m + 1}\times \T^m~:~\tau =T\Big\},\label{himsf:jumpSet}
\end{align}
\end{subequations}
where $x_2\coloneqq (\hat{u}, p, \tau, \mu)$, $k\in\R_{>0}$ is a tunable gain. This controller resets the states $p$ and $\tau$ via \eqref{himsf:jumpMap} every time the timer $\tau$ satisfies $\tau=T$. The constants $0<T_0<T$ are tunable parameters that characterize the frequency of the resets. The parameter $\sigma\in\{0,1\}$ describes the resetting policy for the state $p$. Namely, when $\sigma=1$, we have that $p^+=p$, while $\sigma=0$ leads to $p^+=q$. 

In contrast to \eqref{gisf}, as $(\varepsilon_p,\varepsilon_a,\varepsilon_0)\to0^+$ the hybrid ISC \eqref{hmisf} will emulate the behavior of a regularized version of Nesterov's accelerated ODE with momentum \cite{ODE_Nesterov}, given by $\ddot{u}+\frac{3}{t}\dot{u}+\nabla\tilde{\varphi}(u)=0$, which achieves rates of convergence of order $\mathcal{O}(1/t^2)$ in convex functions, or $\mathcal{O}(e^{-\sqrt{\kappa}t})$ with suitable resets corresponding to $\sigma=0$ in \eqref{himsf:jumpMap}. These resets are similar in spirit to ``restarting'' techniques used in the literature of machine learning \cite{Candes_Restarting}. In the context of model-free feedback control, the resets guarantee enough regularity and robustness in the controller so that it can be interconnected with a dynamical plant in the loop \cite{zero_order_poveda_Lina}. Thus, the hybrid controller is also able to achieve incentive seeking.
\begin{theorem}\label{thm:hmisf}
Suppose that Assumptions \ref{highway:vc}, \ref{assump:performanceIndex} (or \ref{assumption:strongConvex}) and \ref{assumption:ditherfreq} hold. are satisfied. Then, the closed-loop system corresponding to Figure \ref{fig:feedback} with ISC given by \eqref{hmisf}, renders PAS the set $\mathcal{A}_2 \coloneqq \set{(\theta, \hat{u}, p, \tau)~:~(\theta, \hat{u})\in \mathcal{A}_0,~p=\hat{u},~\tau\in [T_0, T]}\times \mathbb{T}^m$ as $(\varepsilon_0,\varepsilon_p, \varepsilon_a)\to 0^+$.\QEDB
\end{theorem}\vspace{-0.4cm}
\begin{proof}
We prove Theorem \ref{thm:hmisf} following a similar approach as in Theorem \ref{thm:gisf}. In particular, first note that the hybrid dynamics \eqref{hmisf} are well-posed in the sense of \cite[Sec. 6]{bookHDS} because the sets $C_2$ and $D_2$ are closed, and the maps $F_2$ and $G_2$ are continuous on these sets. Moreover, neglecting the socio-technical dynamics, and using Assumption \ref{assumption:ditherfreq}, the average dynamics of \eqref{hmisf} correspond to an $\mathcal{O}(\varepsilon_a)$-perturbed version of the hybrid Nesterov gradient dynamics studied in \cite{zero_order_poveda_Lina} for the model-free optimization of static maps. Under Assumptions \ref{highway:vc}, \ref{assump:performanceIndex} or \ref{assumption:strongConvex}, and \ref{assumption:ditherfreq}, these average hybrid dynamics render the set $\{(\hat{u},p,\tau):~\hat{u}=p=\text{argmin}~\tilde{\varphi}(u),\tau\in[T_0,T]\}$ (locally) practically asymptotically stable. By using, sequentially, averaging and singular perturbation theory for perturbed hybrid systems \cite[Thm. 7]{zero_order_poveda_Lina}, we obtain the desired result for the interconnection between the controller and the socio-technical dynamics, which are stable under Assumption \ref{highway:vc}. 
\end{proof}

The key advantage of the ISC \eqref{hmisf} is the incorporation of \emph{dynamic momentum} via the states $(p,\tau)$, as well as periodic resets with frequency dependent on the pair $(T_0,T)$. Note that ``optimal'' restarting frequencies can be used as in \cite{zero_order_poveda_Lina} to avoid oscillations in the control action induced by the presence of momentum. It is well-known that momentum-based optimization algorithms can significantly improve the transient performance in problems where the cost $\tilde{\varphi}$ exhibits shallow convexity properties (e.g., $\kappa\ll 1$). As shown later in Section \ref{sec:numericalExamples}, this will be the case under certain operation conditions of the highway networks.  

On the other hand, when the steady state performance function $\tilde{\varphi}$ is strongly convex and its curvature is not necessarily small, one might wonder if it is possible to achieve better transient performance using non-smooth re-scaled gradient-based dynamics. We investigate this scenario in Section 4.3.
%
%

%

\subsection{Fixed-Time Incentive Seeking Control}
\begin{figure}[t]
    \centering
    \includegraphics[width=0.8\linewidth]{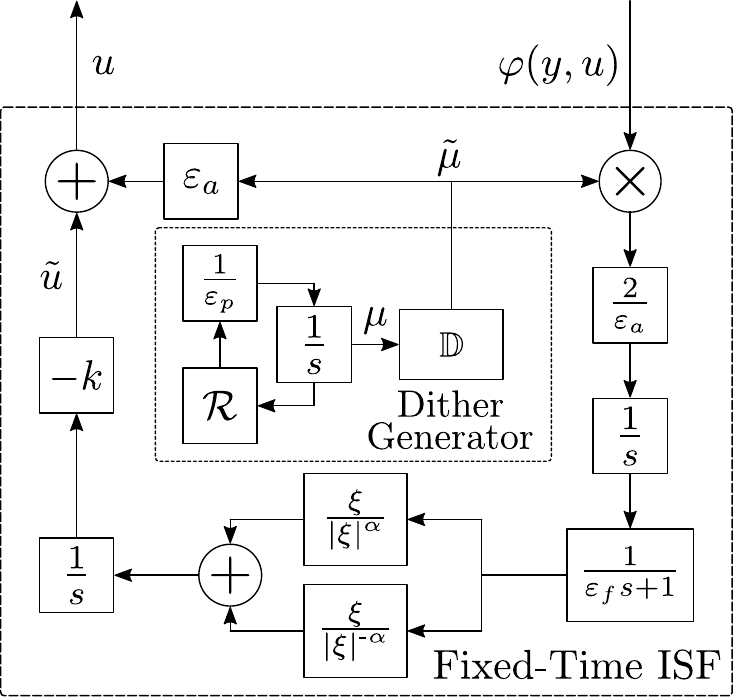}\vspace{-5pt}
    \caption{Fixed-Time Incentive Seeking Control}
    \label{fig:fxes_scheme}
\end{figure}
We now consider the \emph{fixed-time ISC} (FxISC) presented in Figure \ref{fig:fxes_scheme} and described by the following dynamics:
\begin{equation}\label{fxisf}
\left(\begin{array}{c}
         \dot{\hat{u}} \\
         \dot{\xi}\\
         \dot{\mu}
    \end{array}\right)= F_3(x_3) \coloneqq \begin{pmatrix}
    -k\left(\dfrac{\xi}{\abs{\xi}^\alpha} + \dfrac{\xi}{\abs{\xi}^{-\alpha}}\right)\\
    \frac{1}{\varepsilon_f}\left(-\xi +\varphi (y,u)M(\mu)\right)\\
    \frac{1}{\varepsilon_p}\mathcal{R}\mu
    \end{pmatrix}
\end{equation}
where $x_3 \coloneqq (\hat{u}, \xi, \mu)$, and $\alpha\in (0,1)$ is a tunable exponent, and where the right hand side of $\dot{u}$ is defined to be zero whenever $\xi=0$. In this controller, we have incorporated a low-pass filter with state $\xi$ and gain $\varepsilon_f^{-1}$, and the nominal incentive $\hat{u}$ is generated by a combination of sub-linear and super-linear feedback parametrized by the constant $\alpha$. Note that the vector field \eqref{fxisf} is continuous but not Lipschitz continuous at $\xi=0$. The controller is designed to emulate the performance of fixed-time gradient flows \cite{garg2020fixed,PovedaKrsticFXT} as $(\varepsilon_0,\varepsilon_a,\varepsilon_p,\varepsilon_f)\to0^+$. This non-smooth ISC also achieves optimal incentive seeking, but it requires (regional) strong convexity of $\tilde{\varphi}$. 
\begin{theorem}\label{thm:fxisf}
Suppose that Assumptions \ref{highway:vc},  \ref{assumption:strongConvex} and \ref{assumption:ditherfreq} hold. Then, the closed-loop system corresponding to Figure \ref{fig:feedback} with ISC given by \eqref{fxisf}, renders PAS the set $\mathcal{A}_3 \coloneqq \mathcal{A}_0\times\set{0}\times \mathbb{T}^m$ as $(\varepsilon_0,\varepsilon_p, \varepsilon_a, \varepsilon_f)\to 0^+$.\QEDB
\end{theorem}\vspace{-0.4cm}
\begin{proof}
Neglecting the socio-technical dynamics, the average dynamics of \eqref{fxisf} are given by a perturbed version of the fixed-time gradient flows studied in \cite{garg2020fixed}. Under Assumption \ref{assumption:strongConvex}, these dynamics render the optimal incentive fixed-time stable. A direct application of averaging theory for non-smooth systems \cite{PovedaKrsticFXT} allows us to conclude practical (with respect to $\tilde{\Lambda}_u$) fixed-time stability for the ISC interconnected with the socio-technical dynamics \eqref{ode:highway}. 
\end{proof}

In contrast to \eqref{gisf} and \eqref{hmisf}, as  $(\varepsilon_0,\varepsilon_p, \varepsilon_a, \varepsilon_f)\to 0^+$, the nonsmooth ISC \eqref{fxisf}
emulates the behavior of gradient flows able to converge to the optimal incentive before a fixed time $T^*=\frac{\pi}{2k\alpha\kappa}$, where $(\alpha,k)$ are tunable parameters of the controller, and $\kappa$ is given by Assumption \ref{assumption:strongConvex}. Such type of behavior cannot be obtained using smooth (i.e., Lipschitz continuous) ISCs. 


%
\section{NUMERICAL EXAMPLES}
\label{sec:numericalExamples}
In this section, we consider particular realizations of the model introduced in Section \ref{sec:model}, as well as numerical examples of the proposed ISCs.
\subsection{Fast Driver Behavior}\label{sec:numericalExamples:1}
\begin{figure}[t!]
    \centering
    \includegraphics[width=0.9\linewidth]{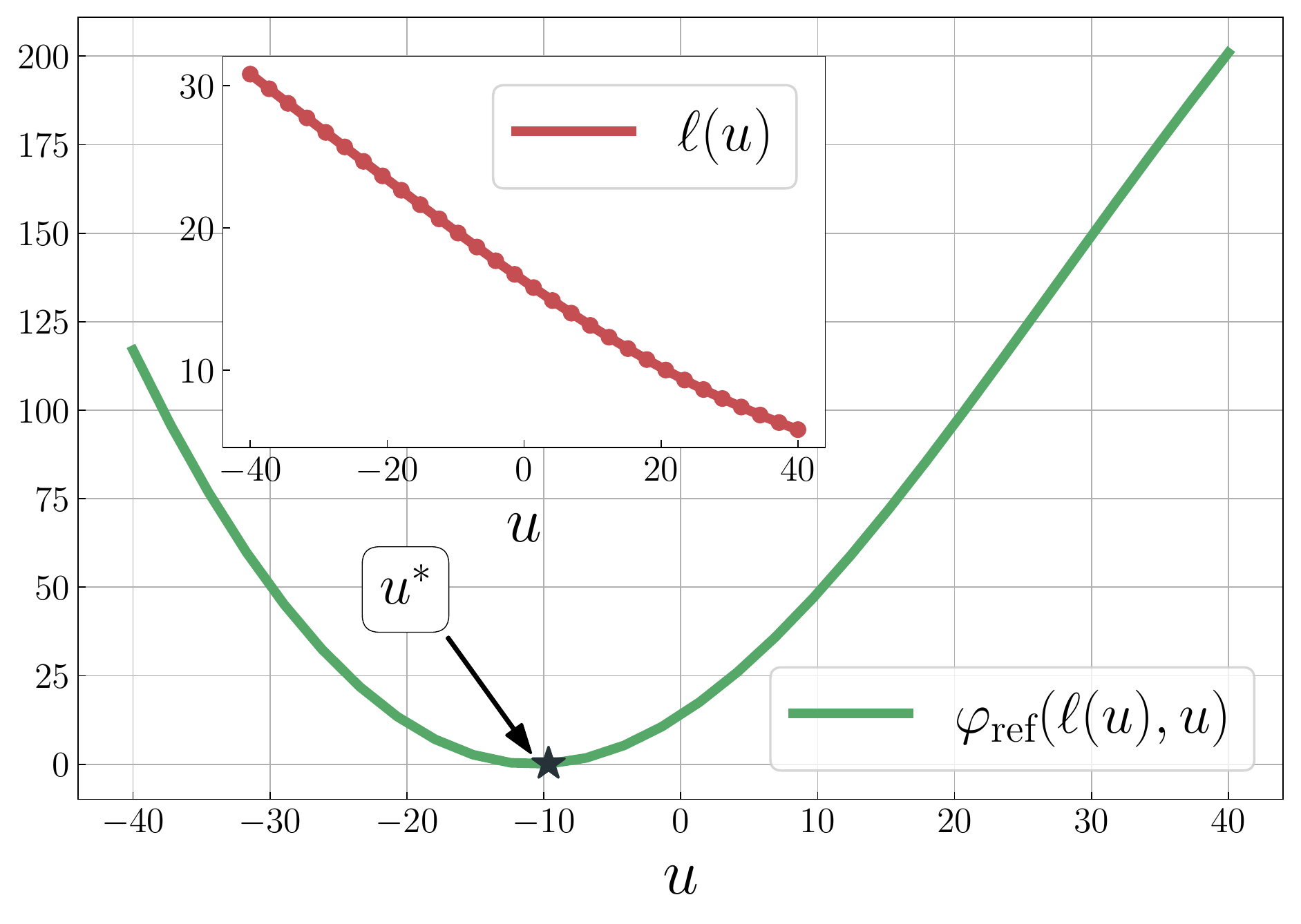}\vspace{-5pt}
    \caption{Suitability of performance index $\varphi_{\text{ref}}$ and response map $\ell(\cdot)$ for the ISC dynamics presented in this paper.}
    \label{fig:viabilityB}
\end{figure}

We first consider a scenario where the driver dynamics are qualitatively faster than the average traffic dynamics. Specifically, we borrow the parameters and structure used in \cite{phan2016model}, based on traffic data of the first dynamic-pricing toll system implemented in the United States: the MnPASS. Thus, the costs of choosing the Express or GP lanes are given by: $c_{\text{EL}}\left(\rho, u\right)= a\frac{L}{\overline{v}(\rho)} + b u + \gamma_{\text{EL}}$, and $c_{\text{GP}}= a\frac{L}{\overline{v}(\rho)}\delta +\gamma_{\text{GP}}$, where $\frac{L}{\overline{v}(\rho)}$ represents the estimated travel time on the Express lane, $\gamma_{\text{EL}},\gamma_{\text{GP}}\in \R$ are offsets used to represent unobservable quantities,  $a,b$ are positive weights, and where $\delta \ge 1$ models the fact that the travel time through the GP lane is assumed to be longer or equal than the one of the Express lane. On the other hand, the macroscopic driver behavior is assumed to have the form: 
\begin{equation}
    q_{\text{EL}}(\rho, u) =\frac{Q}{1 + \text{exp}\left(c_{\text{EL}}(\rho, u)- c_{\text{GP}}(\rho)\right)},\label{qinstatic:numerical}
\end{equation}
which is a logistic function of the marginal cost of choosing the Express lane over the GP lane. The choice of static map in \eqref{qinstatic:numerical} implies that whenever the perceived cost $c_{\text{EL}}$ of choosing the Express lane is lower than the cost $c_{\text{GP}}$ of choosing the GP lane,  the input flow of vehicles to the Express lane will increase. Moreover, we note that an equal inflow of vehicles to the Express and GP lanes is achieved when the marginal cost is equal to zero. To simulate the ISCs we use the parameters $a=0.334,b=0.335,\gamma_{\text{EL}}=1.71781, \gamma_{\text{GP}}=0,v_{\text{jam}}=5\text{[mph]},v_{\text{free}}=65\text{[mph]},\rho_{\text{jam}}=80\text{[veh/mi]}, \rho_{\text{critical}}=25 \text{[veh/mi]}$ and $L=0.7\text{[[mi]}$. To establish a reference density that guarantees free-flow conditions $\overline{v}(\rho_{\text{ref}})\approx v_{\text{free}}$ with a moderate occupation of the lane so that the system is not underutilized, we consider the reference seeking performance index $\varphi_{\text{ref}}$ given in \eqref{reference:seekingcost}, with $\rho_{\text{ref}}= 0.8\rho_{\text{critical}}=20 \text{ [veh/mi]}$. Furthermore, we fix the demand to be $Q = 2170$ vehicles per hour and set $\tilde{\tau}=3$. Using these parameters, we conduct a numerical study that verifies that Assumptions \ref{highway:vc}-\ref{assumption:strongConvex} are satisfied:
%
\begin{figure}[t]
    \centering
    \includegraphics[width=0.95\linewidth]{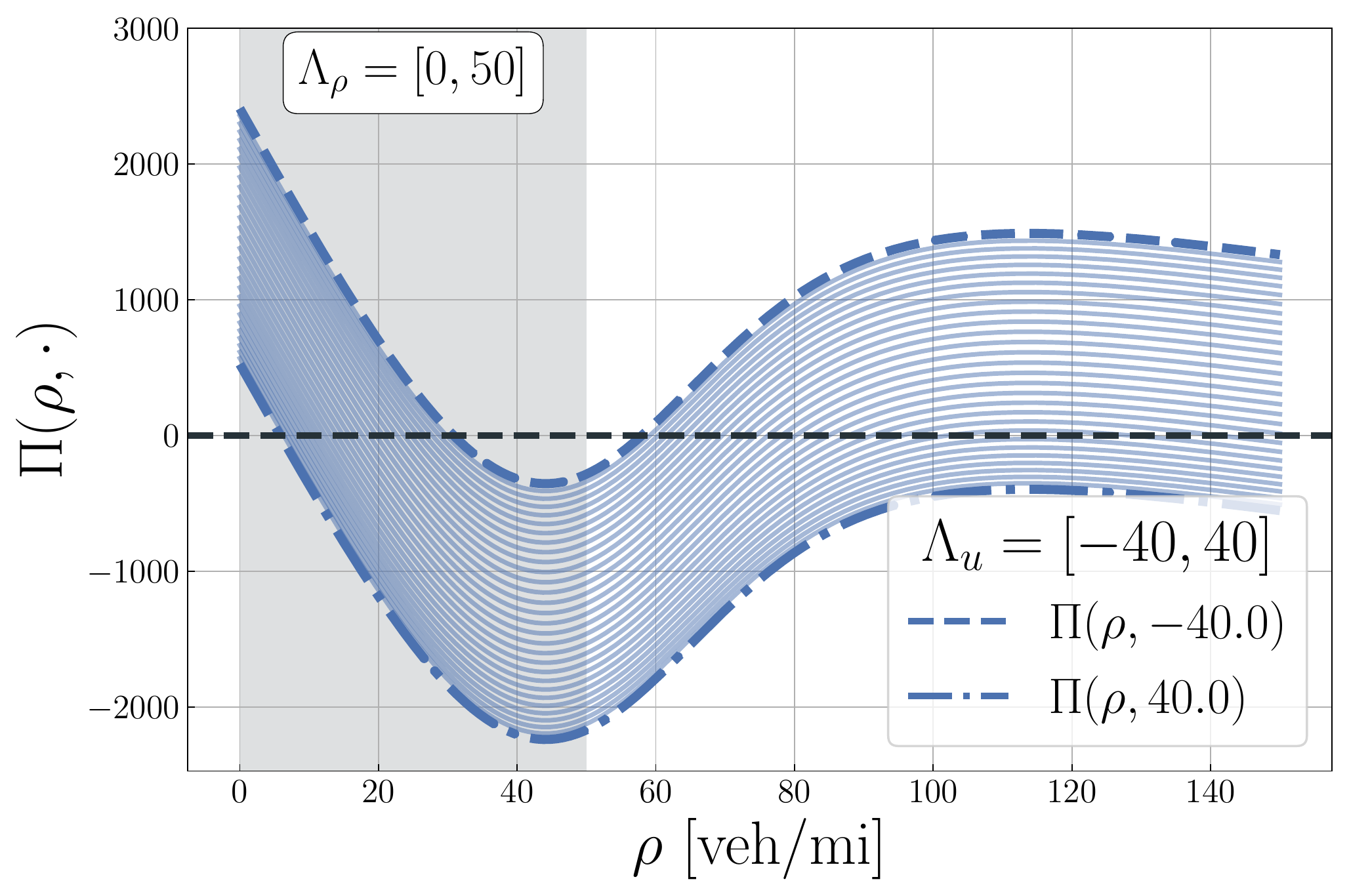}\vspace{-5pt}
    \caption{Evaluation of viability conditions for highway system with fast driver behavior and static input flow map $q_{\text{EL}}$.}
    \label{fig:viabilityA}
\end{figure}
First, we plot in Figure \ref{fig:viabilityA} the vector field \eqref{dynamics:model} for different values of $u\in\Lambda_u\coloneqq [-40, 40]$ and with $\rho=\theta\in [0, 150]$. By arguing graphically, we note that there exists a compact set (interval) $\Lambda_\rho\coloneqq [0, 50]$, such that restricted to values $(\rho,u)\in \Lambda_\rho\times \Lambda_u$, for each $u$ there exists a unique asymptotically stable equilibrium $\rho^*(u)\in \Lambda$. Therefore, we can define the function $\ell:\Lambda_u\to\R$ as $\ell(u)=\rho^*(u)$, which is shown in Figure \ref{fig:viabilityB}. 
Since $\Pi$ in \eqref{dynamics:model:qinstatic} is locally Lipschitz, the previous arguments imply that the socio-technical model of Section \ref{sec:model}, with the particular parameters listed above, satisfy the conditions of Assumption \ref{highway:vc}. On the other hand, $\varphi_{\text{ref}}$ satisfies Assumptions \ref{assump:performanceIndex} and \ref{assumption:strongConvex} by construction and the convexity of $\tilde{\varphi}(u)=\varphi(\ell(u),u)$ in $\Lambda_u$. 

\begin{figure}[t]
    \centering
    \includegraphics[width=0.9\linewidth]{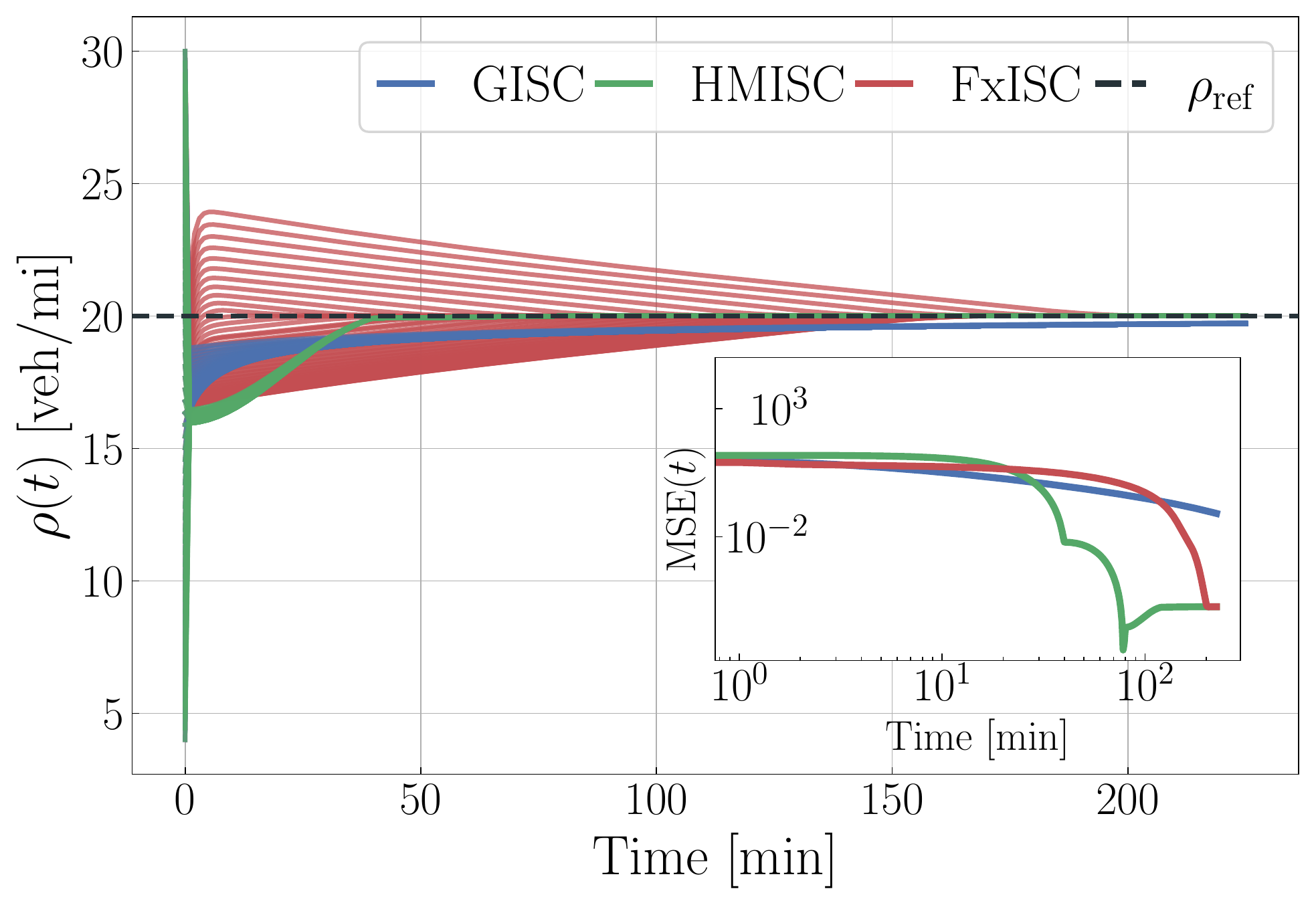}\vspace{-5pt}
    \caption{Trajectories corresponding to 60 different initial conditions of the average density $\rho$ sampled uniformly between $\rho(0)=4$ and  $\rho(0)=30$ vehicles per mile. The inset shows the MSE as a function of time of the 60 trajectories.}
    \label{fig:initialconds}
\end{figure}

\begin{figure}[t]
    \centering
    \includegraphics[width=0.88\linewidth]{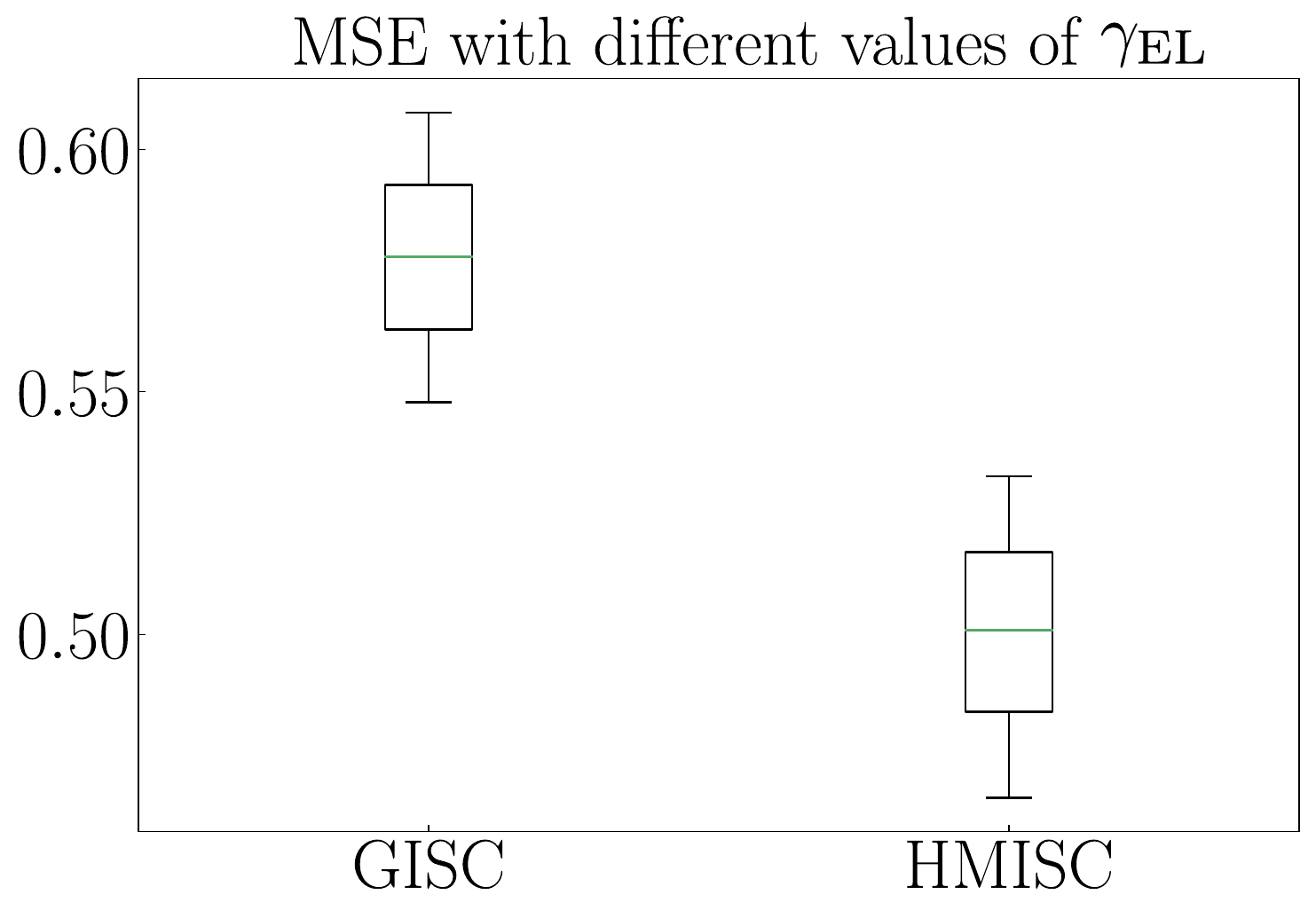}
    \caption{Distribution of the mean of the time-averaged MSE for different values of $\gamma_\text{EL}$ sampled uniformly within $15\%$ of its nominal value.}
    \label{fig:parameterChange}
    \vspace{-0.45cm}
\end{figure}

Following the closed-loop structure of Figure \ref{fig:feedback}, we implement the different ISCs introduced in Section \ref{sec:isf}  interconnected with the highway dynamics of \eqref{dynamics:model:qinstatic}.
For all the controllers, we set $k=1$ and use the dithering frequency $\omega=1$. In the case of the hybrid ISC we choose $\sigma = 0, T_0=0.1$ and $T=20$. For the non-smooth ISC we use $\alpha=0.5$. We simulate the trajectories of the closed-loop systems using $\varepsilon_f = 1,~\varepsilon_a=0.1,$ and $\varepsilon_p = 0.01 $. These parameters guarantee enough time-scale separation between the different elements of the controller, and also between the controller and the highway-dynamics. We uniformly sample $60$ different initial conditions for $\rho$ between $4$ and $30$. We use the initial incentive $u(0)=1$, and we plot the resulting trajectories in Figure \ref{fig:initialconds}. Additionally, we compute and plot the mean squared error (MSE)  $\text{MSE}(t) = \frac{1}{60}\sum_{i=1}^{60} \abs{\rho_i(t)-\rho_{\text{ref}}}^2,$ where $\rho_i$ corresponds to the trajectory resulting from the $i$-th initial condition. As shown in the figure, the hybrid and non-smooth ISCs significantly outperform the smooth ISC algorithm \eqref{gisf}. Note that the hybrid algorithm \eqref{hmisf} generates the typical oscillatory behavior observed in momentum-based algorithms when the damping is sufficiently small. Note also that the hybrid controller seems to generate better transient performance compared to \eqref{fxisf}, since in certain cases it generates a smaller overshoot. The inset of Figure \ref{fig:initialconds} shows that the control signals converge to a small neighborhood of the optimal incentive.  Finally, we note that all the ISCs studied in this paper are well-posed by construction, and therefore they are robust with respect to small bounded additive disturbances acting on the states and dynamics \cite[Thm. 7.21]{bookHDS}. Moreover, their model-free nature allows them to retain their stability and convergence properties when the parameters of the highway change (slowly) over time. For example, Figure \ref{fig:parameterChange} shows the impact of variations on the parameter $\gamma_{\text{EL}}$. Here, we sampled uniformly $20$ different values of $\gamma_{\text{EL}}$ with a maximum variation of $15\%$ with respect to the nominal value $1.71781$, and we simulated the closed-loop dynamics for each one of the ISCs. For each value of $\gamma_{\text{EL}}$ we computed the mean of the time-average MSE, $\overline{\text{tMSE}} \coloneqq \frac{1}{t_f}\int_{0}^{t_f}\text{MSE}(\tau)d\tau$, where $t_f=225\text{[min]}$ is the final time of a simulation run,  over $5$ trajectories obtained by choosing different initial conditions for $\rho$ on the range $[4, 30]$. As seen in Figure \ref{fig:parameterChange}, the results are consistent with the previous results shown in Figure \ref{fig:initialconds}.
%
%
%
%
%
%
\subsection{Dynamics for Macroscopic Driver Behavior}
\begin{figure*}[t]
    \centering
    \includegraphics[width=0.9\linewidth]{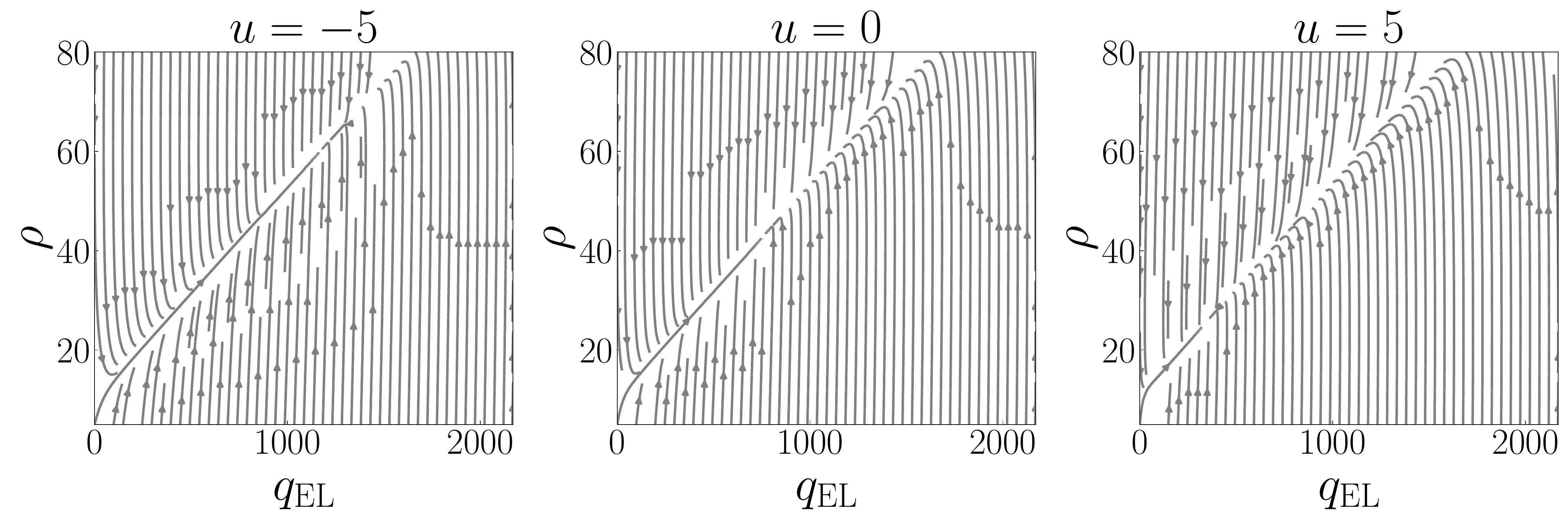}\vspace{-10pt}
    \caption{Evaluation of the viability conditions for highway system with dynamics describing the macroscopic driver behavior, based on the phase plane of the system.}
    \label{fig:viabilityPopDynamics}
\end{figure*}

\begin{figure}[t]
    \centering
    \includegraphics[width=0.85\linewidth]{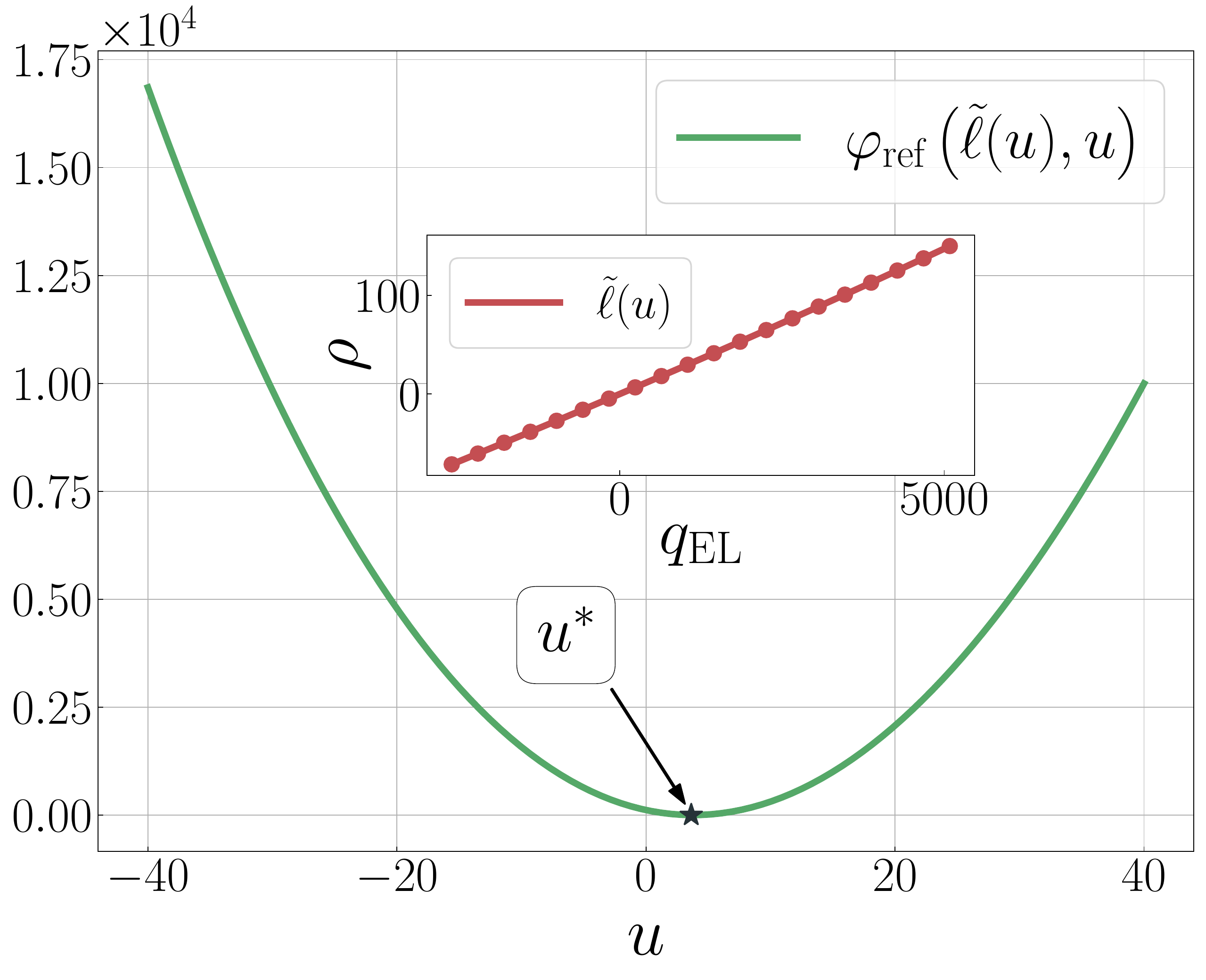}
    \caption{Suitability of performance index $\varphi_{\text{ref}}$ and response map $\tilde{\ell}$ for the ISC with highway model including dynamics for the macroscopic driver behavior. The inset shows the response function $\tilde{\ell}$ projected in the phase plane $\rho~\text{vs.}~ q_{\text{EL}}$ for values of $u\in[-40,40]$.  }
    \label{fig:viabilityPopResponse}
\end{figure}

\begin{figure}[t]
    \centering
    \includegraphics[width=0.95\linewidth]{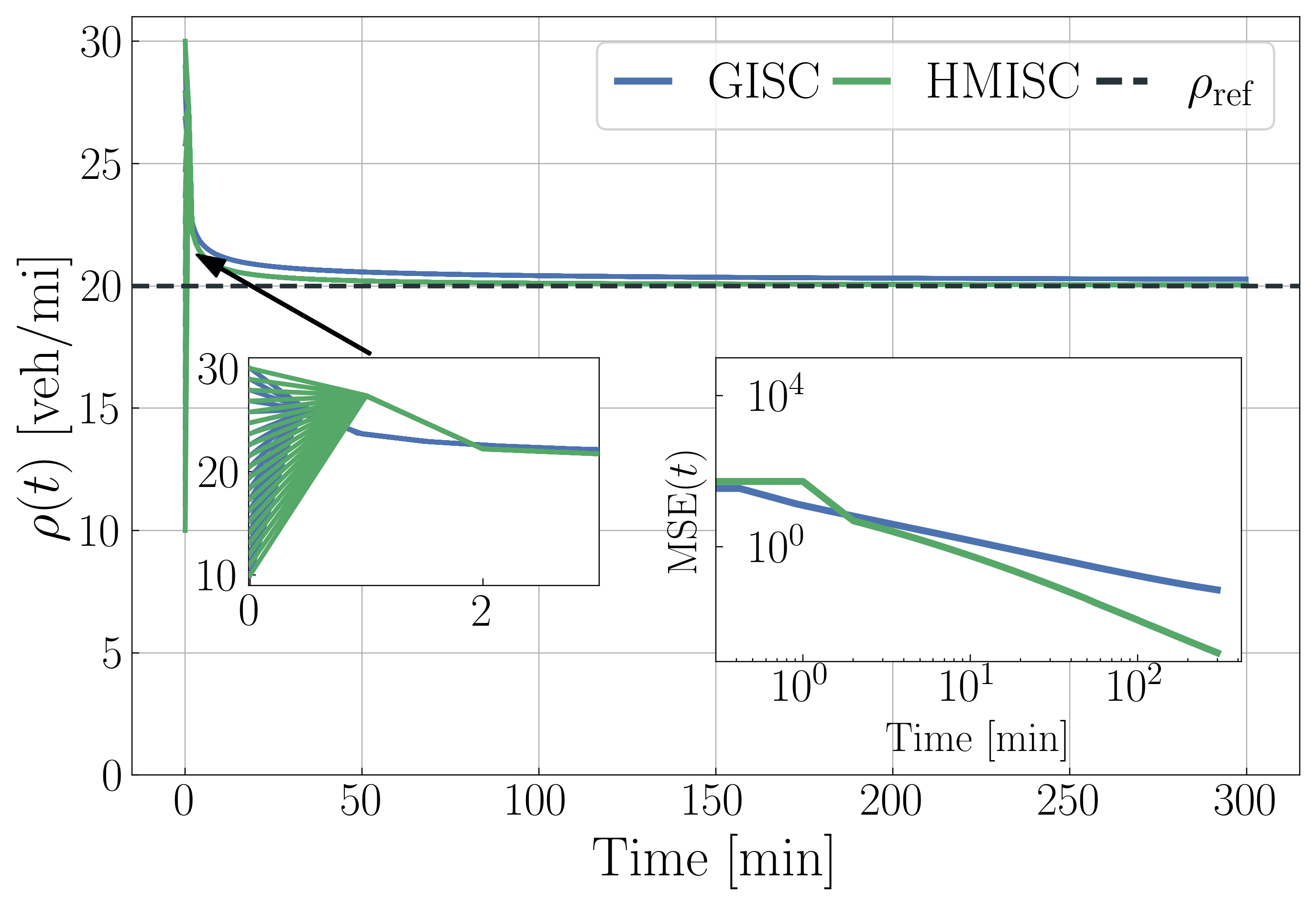}
    \caption{Trajectories resulting from the application of the GISC and HMISC to the highway model with dynamics for the macroscopic driver behavior, and corresponding to 20 different initial conditions of the average density $\rho$ sampled uniformly between $\rho(0)=10$ and  $\rho(0)=30$ vehicles per mile. The inset shows the MSE as a function of time of the 20 trajectories.}
    \label{fig:resultsMacroDynamics}
\end{figure}

When the macroscopic driver behavior and the average vehicle density evolve in similar time scales, the dynamics capturing the driver's response due to changes in the incentives need to be considered in the closed-loop system. In this case, we consider the following dynamics to describe the evolution of $q_{\text{EL}}$:
\begin{equation}\label{logistic:growth:equation}
    \dot{q}_{EL} = \Psi(q_{\text{EL}}, u) \coloneqq -\bigg(c_{\text{EL}}(q_{\text{EL}},u)-c_{\text{GP}}\left(q_{\text{EL}}\right)\bigg),
\end{equation}
meaning that the rate of change of the input flow of vehicles to the Express lane is directly determined by the marginal cost of choosing that lane over the GP lane. Consequently, when the perceived cost $c_{\text{EL}}$ of choosing the Express lane is lower than the cost $c_{\text{GP}}$ of choosing the GP lane, the rate of growth will be instantaneously positive, thus increasing the input flow of vehicles to the Express lane. 
We considered the marginal cost 
\begin{align*}
    c_{\text{EL}}(q_{\text{EL}},u) - c_{\text{GP}}(q_{\text{EL}}) = \left(q_{\text{EL}} -\frac{Q}{2}\right) + \tilde{a}u,
\end{align*}
where $\tilde{a}>0$. For our simulations we use $\tilde{a}=100$ and the same values of $v_{\text{jam}},v_{\text{free}},\rho_{\text{jam}}, \rho_{\text{critical}},Q$ and $\rho_{\text{ref}}$ considered in Section \ref{sec:numericalExamples:1}.  To study the stability properties of the equilibrium points of the highway dynamics, we analyze the phase planes associated to system \eqref{dynamics:model} using $\Phi = \Psi$, $k_m = 1, k_\rho = 1$, and $u\in\Lambda_u=[-40,40]$, and shown in Figure \ref{fig:viabilityPopDynamics} for three particular values of $u$. In all cases, there exists a compact set $\Lambda_{\theta} \subset [0,160]\times[0,Q]$, such that for each $u$ there exists a unique asymptotically stable equilibrium $\theta^*(u)\in \Lambda_{\theta}$. The same property was numerically confirmed to hold for every $u\in\Lambda_u$ by studying the phase plane plots associated to equally spaced inputs taken from $\Lambda_u$, and using continuity of \eqref{dynamics:model}. Therefore, we can define a response function $\tilde{\ell}:\Lambda_u\to \R^2$ by letting $\tilde{\ell}(u)=\theta^*(u)$. The remaining conditions of Assumption \ref{highway:vc}, and Assumptions \ref{assump:performanceIndex} and \ref{assumption:strongConvex} are verified to hold by following analogous graphical arguments to the ones described in Section \ref{sec:numericalExamples:1}. Indeed, Figure \ref{fig:viabilityPopResponse} shows the corresponding plots describing the response function $\tilde{\ell}$ and the performance index  $\varphi_{\text{ref}}\left(\tilde{\ell}(u),u\right)$. In this case, we focus our attention on the two ISCs that showed the best performance in Section \ref{sec:numericalExamples:1}, and
implement the closed-loop structure of Figure \ref{fig:feedback} using the GISC and HMISC. For both controllers we set $\varepsilon_0 = 0.1, 
\varepsilon_\mu = 0.01, \varepsilon_a = 0.001$, $k =  0.01$, and $\omega=1$. For the HMISC, we chose $\sigma = 1, T_0 =0.01$ and $T=0.5$. The control parameters are selected to guarantee enough time-scale separation between the different elements of the controller and the highway-dynamics. As seen in Figure \ref{fig:resultsMacroDynamics}, where we plotted the MSE corresponding to $20$ different trajectories satisfying $u(0)=1$ and $q_{\text{EL}}(0)=Q/3$, the HMISC outperforms the GISC, although in a mildly less dominant fashion to what was observed in Section \ref{sec:numericalExamples:1}. This reduction in the performance gap between the GISC and the HMISC, can be mainly attributed to the fact that, in this setup, the strong convexity parameter is relatively high, meaning that the transient performance increase attained via momentum-based dynamics is not as evident as when the costs exhibit shallow convexity properties.

Finally, we note that in all our numerical experiments the ISCs were tuned to guarantee that the drivers have enough time to react to changes on the incentives induced by the exploratory signal $\tilde{\mu}$ used by the controllers. This behavior is needed to guarantee real-time learning via feedback measurements of the output of the highway network, and it has also been studied in algorithms based on adaptive control \cite{poveda2017class} and reinforcement learning \cite{zhu2015reinforcement}, to name just a few. Potential extensions that could relax these real-time exploration requirements might be studied in the future by incorporating historical data into the controllers, which can be periodically updated during days or weeks to retain sub-optimality of the incentives. Such controllers will naturally be modeled as hybrid dynamical systems.   
%
\section{CONCLUSIONS}
We introduced a new class of incentive-seeking controllers (ISCs) that can learn optimal incentives using only output measurements from traffic in transportation systems, while simultaneously guaranteeing closed-loop stability. We illustrated the benefits of the proposed controllers via numerical experiments in a socio-technical model of a highway system with managed lanes, including the advantages of using nonsmooth and hybrid controllers. The algorithms are agnostic to the exact model of the highway, and robust to small additive disturbances. Future research directions will focus on incorporating past recorded data to minimize real-time exploration in the controllers.


%
\bibliographystyle{ieeetr}
\bibliography{references}

\end{document}